\documentclass{Poitiers}
\usepackage{amsmath}
\usepackage{epsfig} 
\usepackage{graphics} 

\def\ppages{1--10}

\newtheorem{theorem}{Theorem}

\theoremstyle{definition}
\newtheorem{definition}{Definition}

\renewcommand{\figurename}{\footnotesize {\scshape Figure}}

\title[Lie group stability of finite difference schemes]
      {Lie group stability of finite difference schemes}

\author[Emma Hoarau, Claire David, Pierre Sagaut and Thi\text{\^ e}n-Hi\text{\^ e}p L\text{\^ e}]{}

\email{emma.hoarau@onera.fr}

\begin{document}
\maketitle

\centerline{\scshape Emma Hoarau \footnotemark[1], Claire David
\footnotemark[2], Pierre Sagaut \footnotemark[2] and Thi\text{\^
e}n-Hi\text{\^ e}p L\text{\^ e} \footnotemark[1]}
\medskip
{\footnotesize \centerline{\footnotemark[1] ONERA, Computational
Fluid Dynamics and Aeroacoustics Department (DSNA)}
  \centerline{BP 72, 29 avenue de la Division Leclerc}
   \centerline{92322 Ch\text{\^ a}tillon Cedex, France}
\centerline{\footnotemark[2] Universit\'e Pierre et Marie
Curie-Paris 6}
  \centerline{Laboratoire de Mod\'elisation en M\'ecanique, UMR CNRS 7607}
   \centerline{Bo\^ite courrier $n^0 162$, 4 place Jussieu, 75252 Paris, cedex 05,
France}

} 

\medskip

\begin{quote}{\normalfont\fontsize{8}{10}\selectfont
{\bfseries Abstract.} Differential equations arising in fluid
mechanics are usually derived from the intrinsic properties of
mechanical systems, in the form of conservation laws, and bear
symmetries, which are not generally preserved by a finite difference
approximation, and leading to inaccurate numerical results. This
paper develops a method that enables us to build a scheme that
preserves those symmetries. The method is based on the concept of
the differential approximation. A comparison of numerical
performance of the invariant schemes, standard ones and higher order
one has been realised for the Burgers equation.
\par}
\end{quote}

\section{Introduction}
Lie groups were introduced by Sophus Lie in 1870 in order to study
the symmetries of differential equations, yielding thus analytical
solutions. Literature provides substantial works and applications,
\cite{Olverappl}, \cite{Ibragimov}. Symmetry groups can be
determined by an automatic procedure, but it often turn out to be
tedious and induce errors. A large amount of packages using symbolic
manipulations of mathematical expressions have been written. We
mention here some of those works: Schwartz \cite{Schwarz}, Vu and
Carminati\cite{VuCarminati}, Herod \cite{Herod}, Baumann
\cite{Baumann}, Cantwell \cite{Cantwell}.

\noindent In this paper we are interesting in the application of the theory
of Lie group to numerical analysis.

\noindent Finite difference equations used to approximate the
solutions of a differential equation generally do not respect the
symmetries of the original equation, and can lead to inaccurate
numerical results. Various techniques, that enable us to build a
scheme preserving the symmetries of the original differential
equation, have been studied. One of these techniques consists in
constructing  an invariant scheme from a given one by applying the
method of the moving frame in \cite{Olverinv}, \cite{Kim}. Another
one consists in constructing an invariant scheme with the help of
the discret invariants of its symmetry group \cite{Dorodnitsyn},
\cite{Dorodnitsynal1}, \cite{Dorodnitsynal2}, \cite{Bakirova},
\cite{Valiquette} and provides the building of symmetry-adapted
meshes, in preserving the differential equation symmetries. This
technique is based on a direct study of the symmetries of difference
equations and lattices.

\noindent Yanenko \cite{Yanenko} and Shokin \cite{Shokin}, proposed
to apply the Lie group theory to finite difference equations by
means of the differential approximation. Thus, they have set down
conditions under which the differential representation of a finite
difference scheme preserves the group of continuous symmetries of
the original differential equation. They provide a dissipative
scheme, which is called invariant scheme. The resulting scheme is
independent of any change of the reference frame, and its
differential representation is invariant under the symmetries of the
original equation. Ames, Postell and Adams \cite{Ames} have already
used the approach of Yanenko and Shokin to present invariant schemes
in which terms are added to the original difference scheme. They
showed that, in specific cases, the invariant scheme is as accurate
as high order numerical methods.

\noindent In this paper, we focus on the last approach. The method is implemented on
some standard schemes for the Burgers equation. A comparison is made between
the numerical solutions of these schemes and the invariant scheme.

\noindent The paper is organized as follows. Definitions and invariance
condition for differential equations are provided in section 2.
Section 3 recalls the approach of Yanenko and Shokin. Section 4 concentrates
on classical schemes. In section 5, we present a method that enables
us to build an invariant scheme with respect to an otherwise lost
symmetry.

\section{Definitions and invariance condition for differential equations}
\noindent A $r$-parameter Lie group $G_r$ of point transformations in the
Euclidean space $\mathcal{E}(x,u)$ can be written under the form:
\footnotesize{
\begin{equation}
\displaystyle G_r=\{x_i^{*}=\phi_i(x,u,a);\
u_j^{*}=\varphi_j(x,u,a),\ i=1,\dots,m;\ j=1,\dots,n\}
\label{eqn:ga}
\end{equation}
} \normalsize
\noindent Consider a system of $l^{th}$-order differential equations:
\footnotesize{
\begin{equation}
\displaystyle \mathcal
F^{\lambda}\big(x,u,u^{(k_1)},u^{(k_1,k_2)},\dots,u^{(k_1\dots
k_l)}\big)=0,\ \ \lambda=1,\dots,q \label{eqn:ED}
\end{equation}
} \normalsize
\noindent Denote by $u^{(k_1\dots k_p)}$ the vector, the components of which are partial derivatives of order
$p$, namely, $u^{(k_1\dots k_p)}_j=\frac{\partial^p u_j}{\partial x_{k_1}\dots\partial x_{k_p}}$ $j=1,\dots,n$
and $k_1,\dots, k_p \in \{1,\dots,m\}$.

\noindent Denote by $x=(x_1,\dots,x_m)$ the independent variables, $u=(u_1,\dots,u_n)$ the dependent variables,
and $(x_{k_1}\dots x_{k_p})$ a set of elements of the independent variables.

\noindent Equation (\ref{eqn:ED}) is a subset of the Euclidean space
$\mathcal{E}\big(x,u,{u^{(k_1)}},\dots,{u^{(k_1\dots k_l)}}\big)$.
In order to take into account the derivative terms involved in the
differential equation, the action of the group $G_r$ of
transformations in the space $\mathcal{E}\big(x,u)$ needs to be
extended to the space of the derivatives of the dependent variables.

\noindent Denote by $\widetilde G^{(l)}_r$ a $r$-parameter Lie group of point
transformation in the space
$\mathcal{E}\big(x,u,{u^{(k_1)}},\dots,{u^{(k_1\dots k_l)}}\big)$ of
the independent variables, dependent variables and the derivative of
the dependent variables with respect to the independent ones.

\noindent The $l^{th}$-prolongation operator of $G_r$ is: \footnotesize{
\begin{eqnarray}
\displaystyle\widetilde{\mathbf{L}}_\alpha^{(l)}=\xi^\alpha_i(x,u)\frac{\partial
}{\partial x_i}+\eta^\alpha_j(x,u) \frac{\partial }{\partial u_j}+\sigma^{\alpha,(k_1)}_j
\frac{\partial }{\partial
{u_j}^{(k_1)}}+\dots+\sigma^{\alpha,(k_1\dots k_l)}_j \frac{\partial
}{\partial {u_j}^{(k_1\dots k_l)}},\\
\displaystyle \ \ i=1,\dots,m;\ j=1,\dots,n;\ \alpha=1,\dots,r.\nonumber
\end{eqnarray}
} \normalsize $\xi^\alpha_i$, $\eta^\alpha_j$, $\sigma^{\alpha,(k_1)}_j$ and
$\sigma^{\alpha,(k_1\dots k_o)}_j$ are given by: \footnotesize{
\begin{eqnarray}
&&\xi^\alpha_i=
\frac{\partial \phi_i}{\partial a_\alpha}\Big|_{a=0}\nonumber,\ \ \ \ \eta^\alpha_j=\frac{\partial\varphi_j}{\partial a_\alpha}\Big|_{a=0}\nonumber,\ \ \ \ \sigma^{\alpha,(k_1)}_j=
\displaystyle{\frac{\mathcal{D}\eta^\alpha_j}{\mathcal{D} x_{k_1}}
-\sum_{i=1}^m \frac{\partial u_j}{\partial x_i} \frac{\mathcal{D}\xi^\alpha_i}{\mathcal{D}x_{k_1}}} \\
\label{eqn:etas}
&&\sigma^{\alpha,(k_1\dots k_o)}_j =
\displaystyle{\frac{\mathcal{D}\sigma^{\alpha,(k_1\dots
k_{o-1})}_j}{\mathcal{D} x_{k_o}^{\ \ \ \ }}-\sum_{i=1}^m
\frac{\partial^{o} u_j}{\partial x_i \partial x_{k_1}\dots \partial
x_{k_{o-1}}} \frac{\mathcal{D}\xi^\alpha_i}{\mathcal{D}x_{k_o}},\ \
o=2,\dots,l}\nonumber
\end{eqnarray}
} \normalsize where:\footnotesize{$\,\displaystyle
\frac{\mathcal{D}}{\mathcal{D}x_k}=\frac{\partial}{\partial
x_k}+\sum_{j=1}^{n}\frac{\partial u_j}{\partial x_k} \frac{\partial
}{\partial u_j}$} \normalsize

\noindent The system of $l^{th}$-order differential equations is invariant
under the group $\widetilde G^{(l)}_r$ if and only if:
\footnotesize{
\begin{equation}
\displaystyle \widetilde{\mathbf{L}}_\alpha^{(l)}\mathcal F^\lambda \Big|_{\mathcal{F}^{\lambda}=0}=0,\
\ \ \alpha=1,\dots,r;\ \lambda=1,\dots,q \label{eqn:deteq1}
\end{equation}
} \normalsize
\section{Lie group for the differential approximation}
\noindent The finite difference scheme, which approximates the differential
system (\ref{eqn:ED}), can be written as: \footnotesize{
\begin{equation}
\displaystyle{\Lambda^{\lambda}(x,u,h,Tu)=0,\ \ \ \lambda=1,\dots,q}
\label{eqn:scheme}
\end{equation}
} \normalsize \noindent where $h=(h_1,h_2,\dots,h_m)$ denotes the
space step vector, and $T=(T_1,T_2,\dots,T_m)$ the shift-operator
along the axis of the independent variables, defined by:
\footnotesize{
\begin{equation}
T_i[u](x_1,x_2,\dots,x_{i-1},x_i,x_{i+1},\dots,x_m)=u(x_1,x_2,\dots,x_{i-1},x_i+h_i,x_{i+1},\dots,x_m).
\end{equation}
} \normalsize
\begin{definition}
The differential equation: \footnotesize{
\begin{eqnarray}
\displaystyle{\mathcal{P}^{\lambda}\big(x,u,u^{(k_1)},\dots,u^{(k_1\dots k_{l'})}\big)}&=& \displaystyle{\mathcal{F}^{\lambda}\big(x,u,u^{(k_1)},\dots,u^{(k_1\dots k_l)}\big)}\nonumber \\
& & +\displaystyle{\sum_{\beta=1}^s\sum_{i=1}^m (h_i)^{l_\beta} \mathcal{R}^{\lambda}_i(x,u,u^{(k_1)},\dots,u^{(k_1\dots k_{{l'}_{\lambda,i}})})},\nonumber \\
& & \displaystyle{\lambda=1,\dots,q};\ l'= max_{(\lambda,i)}
{l'}_{\lambda,i} \label{eqn:diffapprox}
\end{eqnarray}
} \normalsize is called the $s^{th}$-order differential
approximation of the finite difference scheme (\ref{eqn:scheme}).
In the specific case $s=1$, the above equation is called the first
differential approximation. \label{def:diffapprox}
\end{definition}

\noindent Equation (\ref{eqn:diffapprox}) is obtained from equation (\ref{eqn:scheme}) by applying Taylor series 
 expansion to the components of $T u$ about the point $x=(x_1,\ \dots,\ x_m)$ 
 and truncating the expansion to a given finite order.
\noindent Denote by $G'_r$ a group of transformations in the space
$\mathcal{E}(x,u,h)$: \footnotesize{
\begin{equation}
\displaystyle G'_r=\{x_i^{*}=\phi_i(x,u,a);\
u_j^{*}=\varphi_j(x,u,a);h_i^{*}=\psi_i(x,u,h,a),\ i=1,\dots,m;\
j=1,\dots,n\}
\end{equation}
} \normalsize by $\mathbf{L_\alpha}'$ the basis infinitesimal
operator of $G'_r$: \footnotesize{
\begin{equation}
\displaystyle{\mathbf{L_\alpha}'=
\xi^\alpha_i(x,u)\frac{\partial
}{\partial x_i}+\eta^\alpha_j(x,u) \frac{\partial }{\partial u_j}+\zeta^\alpha_i(x,u,h)\frac{\partial}{\partial
h_i},\ \ \ \alpha=1,\dots,r}
\end{equation}
} \normalsize where \footnotesize{
\begin{equation}
\displaystyle{\zeta^\alpha_i=\frac{\partial\psi_i}{\partial
a_\alpha}\Big{|}_{a=0},\ \ \alpha=1,\dots,r}
\end{equation}
} \normalsize
\noindent and by $\widetilde{G}^{(l')}_r$ a group of transformation in the space $\mathcal{E}(x,u,h,u^{(k_1)},\dots,u^{(k_1\dots k_{l'})})$.\\
The ${l'}^{th}$-prolongation operator of $G'_r$,
$\widetilde{\mathbf{L}}_\alpha^{(l')}$ can be written as: \footnotesize{
\begin{equation}
\displaystyle
\widetilde{\mathbf{L}}_\alpha^{(l')}=\mathbf{L_\alpha}'+\sum_{j=1}^n\sum_{p=1}^{l'}\sigma_j^{\alpha,{(k_1\dots
k_p)}}\frac{\partial }{\partial u_j^{(k_1\dots k_p)}}
\end{equation}
} \normalsize
\begin{theorem}
\noindent The differential approximation (\ref{eqn:diffapprox}) is invariant
under the group $\widetilde G^{(l')}_r$ if and only if \footnotesize{
\begin{equation}
\displaystyle{\widetilde{\mathbf{L}}_\alpha^{(l')}\mathcal{P}^{\lambda}\big((x,u,u^{(k_1)},\dots,u^{(k_1\dots
k_{l'})}\big)\Big|_{\mathcal{P}^{\lambda}=0}=0,\ \ \
\alpha=1,\dots,r;\ \lambda=1,\dots,q}
\end{equation}
} \normalsize or \footnotesize{
\begin{equation}
\displaystyle{\Big[\widetilde{\mathbf{L}}_\alpha^{(l)}\mathcal{F}^{\lambda}+\widetilde{\mathbf{L}}_\alpha^{(l')}\Big(\sum_{\beta=1}^s\sum_{i=1}^m
(h_i)^{l_\beta}
\mathcal{R}^{\lambda}_i\Big)\Big]\Big|_{\mathcal{P}^{\lambda}=0}=0,\
\ \ \alpha=1,\dots,r;\ \lambda=1,\dots,q} \label{eqn:deteq2}
\end{equation}
} \normalsize \label{th:invariance2}
\end{theorem}
\noindent\textbf{Theorem {\ref{th:invariance2}}} provides the equations which enable us to
obtain the symmetry groups of the differential approximation. The unknowns are the infinitesimal
functions $\eta^\alpha_j$, $\xi^\alpha_i$ and $\zeta^\alpha_i$, $i=1,\dots,m$; $j=1,\dots,n$.
The infinitesimals $\sigma^{\alpha,(k_1,\dots,k_o)}_j$, $j=1,\dots,n$, are functions of the partial
derivatives of $\eta^\alpha_j$ and $\xi^\alpha_i$.\\
Equation (\ref{eqn:deteq2}) is simplified by means of the condition
(\ref{eqn:diffapprox}). They lead to an overdetermined system of
differential equations, the unknowns of which are the infinitesimal
functions.
\newpage
\section{The specific case of the Burgers equation}
\subsection{Symmetries of the Burgers equation}
\noindent The Burgers equation can be written as: \footnotesize{
\begin{equation}
\displaystyle{\mathcal{F}(x,t,u,\nu,u_x,u_t,u_{xx})=u_t+u\ u_x-\nu\ u_{xx}=0} \label{eqn:burger}
\end{equation}
} \normalsize \noindent where $\nu\geq 0$ is the dynamic viscosity.

\noindent Denote by $G$ a group of transformations of the Burgers equation in the space $\mathcal{E}(x,t,u,\nu)$
of the independent variables $(x,t)$, the dependent variable $u$, and the viscosity $\nu$. The viscosity is
taken as a symmetry variable in order to enable us to take into account variations of the Reynolds number.\\
$G$ is a set of transformations acting smoothly on the space $\mathcal{E}(x,t,u,\nu)$. 

\noindent The six-dimensional Lie algebra of the group $G$ is generated by the
following operators: \footnotesize{
\begin{eqnarray}
& &\displaystyle{\mathbf{L}_1=\frac{\partial}{\partial x}},\ \displaystyle{\mathbf{L}_2=\frac{\partial}{\partial t}},\ \displaystyle{\mathbf{L}_3=x\frac{\partial}{\partial x}+2 t\frac{\partial}{\partial t}-u\frac{\partial}{\partial u}}\nonumber\\
& &\displaystyle{\mathbf{L}_4=x t\frac{\partial}{\partial x}+t^2\frac{\partial}{\partial t}+(-u t+x)\frac{\partial}{\partial u}},\
\displaystyle{\mathbf{L}_5=t\frac{\partial}{\partial x}+\frac{\partial}{\partial u}},\
\displaystyle{\mathbf{L}_6=-t\frac{\partial}{\partial t}+u\frac{\partial}{\partial u}+\nu\frac{\partial}{\partial \nu}}\label{eqn:op-burgers}
\end{eqnarray}
} \normalsize which respectively correspond to:
\begin{itemize}
\item the space translation
\footnotesize{$:\,(x,t,u,\nu) \longmapsto
(x+\epsilon_1,t,u,\nu)$}\normalsize;
\item the time translation
\footnotesize{$:\,(x,t,u,\nu) \longmapsto
(x,t+\epsilon_2,u,\nu)$}\normalsize;
\item the dilatation
\footnotesize{$:\,(x,t,u,\nu) \longmapsto (\epsilon_3 x,\epsilon_3^2 t,\epsilon_3^{-1}
u,\nu)$}\normalsize;
\item the projective transformation
\footnotesize{$:\,(x,t,u,\nu) \longmapsto \Big(\frac{x}{1-\epsilon_4 t},\frac{t}{1-\epsilon_4 t},x\epsilon_4+u(1-\epsilon_4 t),\nu\Big)$}\normalsize;
\item the Galilean transformation
\footnotesize{$:\,(x,t,u,\nu) \longmapsto (x+\epsilon_5\ t,t,u+\epsilon_5,\nu)$}\normalsize;
\item the dilatation
\footnotesize{$:\,(x,t,u,\nu) \longmapsto (x,\epsilon_6^{-1} t,\epsilon_6 u, \epsilon_6\nu)$}\normalsize.
\end{itemize}
\normalsize
{$(\epsilon_i)_{i=1,\dots,6}$} \normalsize are constants.
\normalsize
\subsection{Symmetries of first differential approximations}
\label{Sym-FDA}
\noindent Denote by $h$ the mesh size, $\tau$ the time step, $N_x$ the number
of mesh points, $N_t$ the number of time steps, and $u^n_i,\ i \in
\{0,\dots,N_t\},\ n \in \{0,\dots,N_x\}$ the
discrete approximation of $u(i h,n \tau)$.

\noindent In order to shorten the size of the finite difference scheme expressions, we use the following notations
introduced by Hildebrand in \cite{Hildebrand}:
\footnotesize{
\begin{eqnarray*}
\displaystyle \delta (u^n_i)=\frac{u^n_{i+\frac{1}{2}}-u^n_{i-\frac{1}{2}}}{h},& & \mu (u^n_i)=\frac{u^n_{i+\frac{1}{2}}+u^n_{i-\frac{1}{2}}}{2}\\
\displaystyle \delta^{+} (u^n_i)=\frac{u^n_{i+1}-u^n_{i}}{h},& & \delta^{-} (u^n_i)=\frac{u^n_{i}-u^n_{i-1}}{h},\,\,\,\,\,\,\,\, E^{\alpha}u^n_{i}=u^n_{i+\alpha}
\end{eqnarray*}
} \normalsize
\noindent The Burgers equation can be discretized by means of:
\begin{itemize}
\item \textbf{the FTCS (forward-time and centered-space) scheme}:
\footnotesize{
$$\displaystyle {\frac{u^{n+1}_i -u^n_i}{\tau}+\frac{\mu\delta}{h}\big(\frac{u^2}{2}\big)^n_{i}-\nu \frac{\delta^2}{h^2}u^n_{i}=0}$$
} \normalsize
\item \textbf{the Lax-Wendroff scheme}:
\footnotesize{
$$\displaystyle {\frac{u^{n+1}_i -u^n_i}{\tau}+\frac{\mu\delta}{h}\big(\frac{u^2}{2}\big)^n_{i}-\nu \frac{\delta^2}{h^2}u^n_{i}+A^n_i=0}$$
} \normalsize
\noindent where:
\footnotesize{
\begin{eqnarray*}
\displaystyle A^n_i=&-&\frac{\tau}{2 h^2}\Big[E^{\frac{1}{2}}u^n_{i}\ \delta^{+}\big(\frac{u^2}{2}\big)^n_{i}-E^{-\frac{1}{2}}u^n_{i}\
\delta^{-}\big(\frac{u^2}{2}\big)^n_{i}\Big]-\frac{\nu^2\tau}{2}\Big[\frac{\delta^4}{h^4}u^n_{i}\Big]\\
&+&\frac{\nu\tau}{2 h^3}\Big[E^{\frac{1}{2}}u^n_{i}\ \delta^{2}(E^{\frac{1}{2}}u^n_{i})-E^{-\frac{1}{2}}u^n_{i}\ \delta^{2}(E^{-\frac{1}{2}}u^n_{i})\Big]
+\frac{\nu\tau}{2}\Big[\frac{\mu \delta^3}{h^3}\big(\frac{u^2}{2}\big)^n_{i}\Big]
\end{eqnarray*}
} \normalsize
\item \textbf{the Crank-Nicolson scheme}:
\footnotesize{
\begin{eqnarray*}
\displaystyle \frac{u^{n+1}_i-u^{n}_i}{\tau}+\frac{\mu\delta}{h}\Big[\big(\frac{u^2}{2}\big)^{n+1}_{i}+\big(\frac{u^2}{2}\big)^n_{i}\Big]-\nu \frac{\delta^2}{h^2}[u^{n+1}_{i}+u^n_{i}]=0
\end{eqnarray*}
} \normalsize
\end{itemize}
\normalsize
\noindent Linear stability properties and the related orders of approximation are displayed in Table \ref{Table1} (where $CFL=\frac{a\tau}{h}$, $S=\frac{\nu\tau}{h^2}$ and $S^*=\big(\nu+\frac{a h CFL}{2}\big)\frac{\tau}{h^2}$).\\

\begin{table}[!hbp]
\begin{center}
\begin{tabular}{|c|c|c|}
  \hline   \textbf{Scheme}& \textbf{Stablility condition}& \textbf{Error}\\
   \hline
   FTCS           & \footnotesize{$S\leq \frac{1}{2}$, $CFL\leq 1$}  & \footnotesize{$ \mathcal{O}(\tau,h^2)$}\\
   \hline
   Lax-Wendroff   & \footnotesize{$S^*\leq \frac{1}{2}$, $CFL\leq 1$}& \footnotesize{$\displaystyle \mathcal{O}(\tau^2,h^2)$}\\
   \hline
   Crank-Nicolson & \footnotesize{unconditional stability}           & \footnotesize{$\displaystyle \mathcal{O}(\tau^2,h^2)$}\\
   \hline
\end{tabular}
\end{center}
\caption{Table of finite difference schemes} \label{Table1}
\end{table}
\normalsize

\noindent Consider ${u_{i}}^{n}$ as a function of the time step $\tau$, and of the mesh size $h$, expand it at a given order by means of its Taylor series, and neglect the $o(\tau^{\alpha})$ and $o({h}^{\beta})$ terms, where $\alpha$ and $\beta$ depend on the order of the schemes. This yields the differential representation of the finite difference equation.\\
The following differential representations are obtained:
\begin{itemize}
\item\textbf{for the FTCS scheme:}
\footnotesize{
\begin{eqnarray*}
\displaystyle u_t+\frac{1}{2}(u^2)_x-\nu\
u_{xx}+\frac{\tau}{2}g_2+\frac{h^2}{12}(u^2)_{xxx}-\frac{\nu
h^2}{12}u_{xxxx}=0
\end{eqnarray*}
} \normalsize
\item \textbf{for the Lax-Wendroff scheme:}
\footnotesize{
\begin{eqnarray*}
\displaystyle u_t+\frac{1}{2}(u^2)_x-\nu\
u_{xx}+\frac{\tau^2}{6}g_3+\frac{h^2}{12}(u^2)_{xxx}-\frac{\nu
h^2}{12}u_{xxxx}=0
\end{eqnarray*}
} \normalsize
\item \textbf{for the Crank-Nicolson scheme:}
\footnotesize{
\begin{eqnarray*}
u_t+\frac{1}{2}(u^2)_x-\nu
u_{xx}+\tau^2\big(\frac{g_3}{6}+\frac{1}{4}(g^2_1+ug_2)_x-\frac{\nu}{4}(g_2)_{xx}\big)+\frac{h^2}{12}(u^2)_{xxx}-\frac{\nu h^2}{12}u_{xxxx}=0
\end{eqnarray*}
} \normalsize
\end{itemize}
\noindent where \footnotesize{$g_1=-\big(\frac{u^2}{2}\big)_x+\nu u_{xx}$,
$g_2=\big(-g_1 u\big)_x+\nu\big(g_1\big)_{xx}$, $g_3=\big(-g_2 u
-g^2_1\big)_x+\nu \big(g_2\big)_{xx}$ } \normalsize \vspace{0.4cm}

\noindent Denote by $G'$ the group of transformations of a first differential approximation in the space $\mathcal{E}(x,t,u,h,\tau,\nu)$ of the independent
variables $(x,t)$ and the dependent variable $u$, the step size
variables $(h,\tau)$ and the viscosity $\nu$.

\noindent The ${l'}^{th}$-prolongation of $G'$ can be written as:
\footnotesize{
\begin{eqnarray}
\displaystyle{\widetilde{\mathbf{L}}_\alpha'^{(l')}}&=&\displaystyle{\xi^\alpha_1\frac{\partial}{\partial
x}+\xi^\alpha_2\frac{\partial}{\partial t}+\eta^\alpha
\frac{\partial}{\partial
u}+\sum_{p=1}^{l'}\sigma_j^{\alpha,(k_1\dots k_p)}\frac{\partial
}{\partial u_j^{(k_1\dots
k_p)}}+\zeta^\alpha_1\frac{\partial}{\partial
h}+\zeta^\alpha_2\frac{\partial}{\partial
\tau}+\theta^\alpha\frac{\partial}{\partial \nu}}
\end{eqnarray}
}\normalsize
\noindent where $l'$ has been defined in \textbf{definition {\ref{def:diffapprox}}}.

\noindent \textbf{Theorem \ref{th:invariance2}} enables us to obtain the necessary
and sufficient condition of invariance of the first differential
approximation $\mathcal{P}$: \footnotesize{
\begin{equation}
\displaystyle{\widetilde{\mathbf{L}}_\alpha'^{(l')}
\mathcal{P}\Big{|}_{\mathcal{P}=0}=0} \label{eqn:condinvdiffapprox}
\end{equation}
} \normalsize
\noindent \textbf{Theorem \ref{th:invariance2}} is applied to the differential representations of the above schemes.

\noindent The resolution of the determining equations of each first
differential approximation yields the $4$-parameter group:
\footnotesize{
\begin{eqnarray}
\displaystyle{\xi^\alpha_1=a+b\ x},&\ \ \displaystyle{\xi^\alpha_2=c+(2 b-d)\  t},&\ \ \displaystyle{\eta^\alpha=(-b+d)\ u}\\
\displaystyle{\zeta^\alpha_1=b\ h},&\ \
\displaystyle{\zeta^\alpha_2=(2 b-d)\ \tau},&\ \
\displaystyle{\theta^\alpha=e \nu}\nonumber
\end{eqnarray}
} \normalsize The $4$-dimensional Lie algebra of $G'$ is generated
by: \footnotesize{
\begin{eqnarray}
& &\displaystyle{\mathbf{L}_1=\frac{\partial}{\partial x}},\ \
\displaystyle{\mathbf{L}_2=\frac{\partial}{\partial t}}\nonumber ,\ \ \displaystyle{\mathbf{L'}_3=x\frac{\partial}{\partial x}+2 t\frac{\partial}{\partial t}-u\frac{\partial}{\partial u}+h\frac{\partial}{\partial h}+2 \tau\frac{\partial}{\partial \tau}}\\
& & \displaystyle{\mathbf{L'}_4=-t\frac{\partial}{\partial
t}+u\frac{\partial}{\partial u}-\tau \frac{\partial}{\partial
\tau}+\nu \frac{\partial}{\partial \nu}}
\end{eqnarray}
} \normalsize These operators are respectively related to:
\begin{itemize}
\item the space translation
\footnotesize{$:\,(x,t,u,h,\tau,\nu) \longmapsto (x+\epsilon_1,t,u,h,\tau,\nu)$}\normalsize ;
\item the time translation
\footnotesize{$:\,(x,t,u,h,\tau,\nu) \longmapsto (x,t+\epsilon_2,u,h,\tau,\nu)$}\normalsize ;
\item the dilatation
\footnotesize{$:\,(x,t,u,h,\tau,\nu) \longmapsto (\epsilon_3 x,\epsilon_3^2 t,\epsilon_3^{-1} u,\epsilon_3 h,\epsilon_3^2 \tau,\nu)$}\normalsize ;
\item the dilatation
\footnotesize{$:\,(x,t,u,h,\tau,\nu) \longmapsto (x,\epsilon_4^{-1} t,\epsilon_4 u,h,\epsilon_4^{-1}\tau,\epsilon_4 \nu)$}\normalsize ;
\end{itemize}
\noindent where $(\epsilon_i)_{i=1,\dots,4}$ are constants.\\
\noindent The above finite difference equations are preserved by the space translation, the time translation and both dilatations.

\noindent Approximating the Burgers equation by the above finite difference
equations results in the loss of the projective and Galilean transformations.

\section{The invariant scheme}
\subsection{Invariant scheme construction}
\noindent An invariant scheme is constructed in such a way that the
related differential approximation preserves the symmetries of the
Burgers equation. \noindent We propose to approximate the Burgers
equation by the following finite difference scheme: \footnotesize{
\begin{eqnarray}
\displaystyle\frac{u^{n+1}_i-u^n_i}{\tau}+\frac{1}{h}\big(\mu\delta-\frac{\mu\delta^3}{6}\big)\big(\frac{u^2}{2}\big)^n_{i}-\nu\frac{1}{h^2}\big(\delta^2-\frac{\delta^4}{12}\big)(u^n_i)-\Big(\Omega^n_{i+\frac{1}{2}}
\delta^{+}-\Omega^n_{i-\frac{1}{2}} \delta^{-}\Big)u^n_i=0
\label{eqn:inv}
\end{eqnarray}
} \normalsize \noindent where $\Omega^n_i=\Omega(x_i,t_n,u^n_i)$ is
defined next so that the related differential representation is
preserved by the symmetries of the Burgers equation. \noindent The
scheme has second-order accuracy in space and first-order accuracy
in time. The derivatives $\big(\frac{u^2}{2}\big)_x$ and $u_{xx}$
are approximated by fourth order accuracy difference expressions:
\footnotesize{
\begin{eqnarray}
\displaystyle\big(\frac{\mu\delta}{h}-\frac{\mu\delta^3}{6 h}\big)(u^n_i)=\big(u_x-\frac{h^4}{30}u_{5x}\big)^n_i+\mathcal{O}(h^6),\
\displaystyle\big(\frac{\delta^2}{h^2}-\frac{\delta^4}{12 h^2}\big)(u^n_i)=\big(u_{xx}-\frac{h^4}{90}u_{6x}\big)^n_i+\mathcal{O}(h^6)
\label{eqn:4-order-acc}
\end{eqnarray}
} \normalsize
\noindent The truncation error of the difference scheme (\ref{eqn:inv}) can be written as:
\footnotesize{
\begin{eqnarray}
\displaystyle \epsilon&=&\frac{\tau}{2}u_{tt}-{h^2}\Big(\Omega u_x\Big)_x+\mathcal{O}(\tau^2)+\mathcal{O}(h^4)\nonumber
\end{eqnarray}
} \normalsize
$u_{tt}$ is replaced by an expression involving partial derivatives with respect to $x$, by using the Burgers equation:
\footnotesize{
\begin{eqnarray}
\displaystyle u_{tt}=(u^2 u_x)_x-\nu(u u_{xx})_x-\nu \big(\frac{u^2}{2}\big)_{xxx}+\nu^2 u_{xxxx}
\end{eqnarray}
} \normalsize
\noindent Replacing the previous expression in the truncation error leads to:
\footnotesize{
\begin{eqnarray*}
\displaystyle \epsilon&=&\Big(C u_x\Big)_x-\frac{\nu\tau}{2}\Big(u u_{xx}\Big)_x-\frac{\nu\tau}{2}\big(\frac{u^2}{2}\big)_{xxx}+\frac{\nu^2\tau}{2}u_{xxxx}+\mathcal{O}(\tau^2)+\mathcal{O}(h^4)
\end{eqnarray*}
} \normalsize
\noindent where $\displaystyle C=\frac{\tau}{2}u^2-h^2\Omega$.\\
\noindent It is convenient for the calculation of $C$ that the truncation error is reduced to:
\footnotesize{
\begin{eqnarray*}
\displaystyle \epsilon&=&\Big(C u_x\Big)_x+\mathcal{O}(\tau^2)+\mathcal{O}(h^4)
\end{eqnarray*}
} \normalsize
\noindent The related finite difference scheme is the following first order accuracy in time and second order accuracy in space:
\footnotesize{
\begin{eqnarray}
\displaystyle
\frac{u^{n+1}_i-u^n_i}{\tau}+\frac{1}{h}\big(\mu\delta-\frac{\mu\delta^3}{6}\big)\big(\frac{u^2}{2}\big)^n_{i}-\nu\frac{1}{h^2}\big(\delta^2-\frac{\delta^4}{12}\big)(u^n_i)-\Big(\Omega^n_{i+\frac{1}{2}}
\Delta_1-\Omega^n_{i-\frac{1}{2}} \Delta_{-1}\Big)u^n_i\nonumber\\
\displaystyle
+\frac{\nu\tau}{2}\Big(u^n_{i+\frac{1}{2}}\frac{\mu\delta^2}{h^2}(u^n_{i+\frac{1}{2}})-u^n_{i-\frac{1}{2}}\frac{\mu\delta^2}{h^2}(u^n_{i-\frac{1}{2}})\Big)-\frac{\nu^2\tau}{2}\frac{\delta^4}{h^4}u^n_i+\frac{\nu\tau}{2}\frac{\mu\delta^3}{h^3}\big(\frac{u^2}{2}\big)^n_i=0
\label{eqn:inv2}
\end{eqnarray}
} \normalsize
\noindent and the differential approximation can be written as:
\footnotesize{
\begin{eqnarray}
\displaystyle {\mathcal{P}(x,t,u,\nu,u_x,u_t,u_{xx})=u_t+u\ u_x-\nu\ u_{xx}+(C u_x)_x=0}
\label{eqn:diff-approx-burgers}
\end{eqnarray}
} \normalsize
\noindent The von Neumann stability analysis of scheme (\ref{eqn:inv2}) under a linearized form provides
the following necessary conditions for $S$, $CFL$ and $\Omega_\tau=\Omega\tau$:
\footnotesize{
\begin{eqnarray}
\displaystyle CFL^2-2S-2\Omega_\tau\leq0,\ \ \ \ 0\leq\frac{4 S}{3}-2 S^2+\Omega_\tau\leq\frac{1}{2}
\label{eqn:stab-cond}
\end{eqnarray}
} \normalsize \noindent If $\Omega$ takes is sufficiently close to
zero, these conditions become then sufficient for the linear formulation.
\subsection{Calculation of the artificial viscosity term}
\noindent Here we describe the method for determining the artificial
viscosity term $(C u_x)_x$, which is constructed in such a way that
the differential approximation (\ref{eqn:diff-approx-burgers}) is
preserved by the symmetries of the Burgers equation. $C$ is a
function of the variables $(x,t,u,\tau,h)$, and also depends on the
partial derivatives of $u$ with respect to $x$: $u_x$ and $u_{xx}$.
$C=C(x,t,h,\tau,u,u_x,u_{xx})$.
The necessary and sufficient condition for the differential approximation to be an invariant of the Burgers equation symmetry group is:
\footnotesize{
\begin{eqnarray}
\displaystyle \widetilde{\mathbf{L}}_\alpha^{(2)} (u_t+u\ u_x-\nu\ u_{xx})\Big{|}_{\mathcal{P}=0}+\widetilde{\mathbf{L}}_\alpha^{(3)} ((C u_x)_x)\Big{|}_{\mathcal{P}=0}=0
\label{eqn:calcul-C}
\end{eqnarray}
} \normalsize
\noindent Equation (\ref{eqn:calcul-C}) provides the determining equations of the symmetry group of equation (\ref{eqn:diff-approx-burgers}).
The determining equations involve partial derivatives of the unknown function $C$ and partial derivatives of the infinitesimal
functions of $G^{'}$, which is the symmetry group of the differential representation of the invariant scheme.

\noindent The infinitesimal functions of $G^{'}$ have the following expressions:
\footnotesize{
\begin{eqnarray}
\displaystyle{\xi^\alpha_1=a+b\ x+c\ t+d\ t x},&\ \displaystyle{\xi^\alpha_2=e+d\ t^2+(2 b-f)\  t},&\ \displaystyle{\zeta^\alpha_1=b\ h},\\
\displaystyle{\zeta^\alpha_2=(2 b-f)\ \tau},&\ \displaystyle{\eta^\alpha=c+d\ x+(-b-d\ t+f)\ u},&\ \displaystyle{\theta^\alpha=f \nu}\nonumber
\end{eqnarray}
} \normalsize
\noindent The determining equation with respect to the unknown function $C$ is simplified in using the infinitesimal functions of each subgroup of $G^{'}$.

\noindent The determining equations of each subgroup of $G^{'}$ provides the following linear partial differential equations and the expressions for
$C$:
\begin{itemize}
\item the space translation
\footnotesize{$\frac{\partial}{\partial x} C=0$ $\Rightarrow$ $C=C_1(t,h,\tau,u,u_x,u_{xx})$}\normalsize;
\item the time translation
\footnotesize{$\frac{\partial}{\partial t} C=0$ $\Rightarrow$ $C=C_2(x,h,\tau,u,u_x,u_{xx})$}\normalsize;
\item the dilatation
\footnotesize{$x\frac{\partial}{\partial x} C+2t\frac{\partial}{\partial t}C-u\frac{\partial}{\partial u} C+h\frac{\partial}{\partial h}
C+2\tau\frac{\partial}{\partial \tau}C=0$ $\Rightarrow$ $C=C_3(\frac{t}{x^2},\frac{h}{x},u x,\frac{\tau}{x^2},u_x,u_{xx})$}\normalsize;
\item the projective transformation
\footnotesize{$\frac{\partial}{\partial x} C=0,\, \frac{\partial}{\partial u}C=0,\, \frac{\partial}{\partial u_{xx}}C=0,\, t^2\frac{\partial}{\partial
t}C+2\frac{\partial}{\partial u_{x}}C=0$ $\Rightarrow$ $C=C_4(h,\tau,\frac{2+t u_x}{t})$}\normalsize;
\item the Galilean transformation
\footnotesize{$\frac{\partial}{\partial u} C+t\frac{\partial}{\partial x} C$ $\Rightarrow$ $C=C_5(\frac{u t-x}{t},t,h,\tau,u,u_x,u_{xx})$}\normalsize;
\item the dilatation
\footnotesize{$-t\frac{\partial}{\partial t}C+u\frac{\partial}{\partial u} C+\nu\frac{\partial}{\partial \nu} C
-\tau\frac{\partial}{\partial \tau}C=0$ $\Rightarrow$ $C=\frac{1}{t}C_6(x,h,\frac{\tau}{t})$}\normalsize.
\end{itemize}
\subsection{Numerical application}
\noindent The numerical resolution of the Burgers equation has been
implemented for scheme (\ref{eqn:inv2}), the standard schemes (cf.
section \ref{Sym-FDA}) and a scheme with second-order accuracy in
time and fourth-order accuracy in space, which is obtained from the
invariant scheme when $C=0$. The solutions are calculated in the
reference frame $(F1)$ and in the one $(F2)$ resulting from the
Galilean transformations $(x,t,u,\nu) \longmapsto (x+t,t,u+1,\nu)$.
The artificial viscosity has the following expression:
\footnotesize{
\begin{eqnarray}
\displaystyle C=-0.01 t(t u-x)^2 (u_x)^2
\end{eqnarray}
} \normalsize \noindent $\displaystyle \Omega=\frac{1}{h^2}(\frac{\tau}{2}u^2-C)$ is in a
sufficiently small neighborhood of zero that we have the sufficiency
of conditions (\ref{eqn:stab-cond}) for the linear formulation.

\begin{figure}[!hbp]
\resizebox*{0.49\columnwidth}{0.23\textheight}{\includegraphics{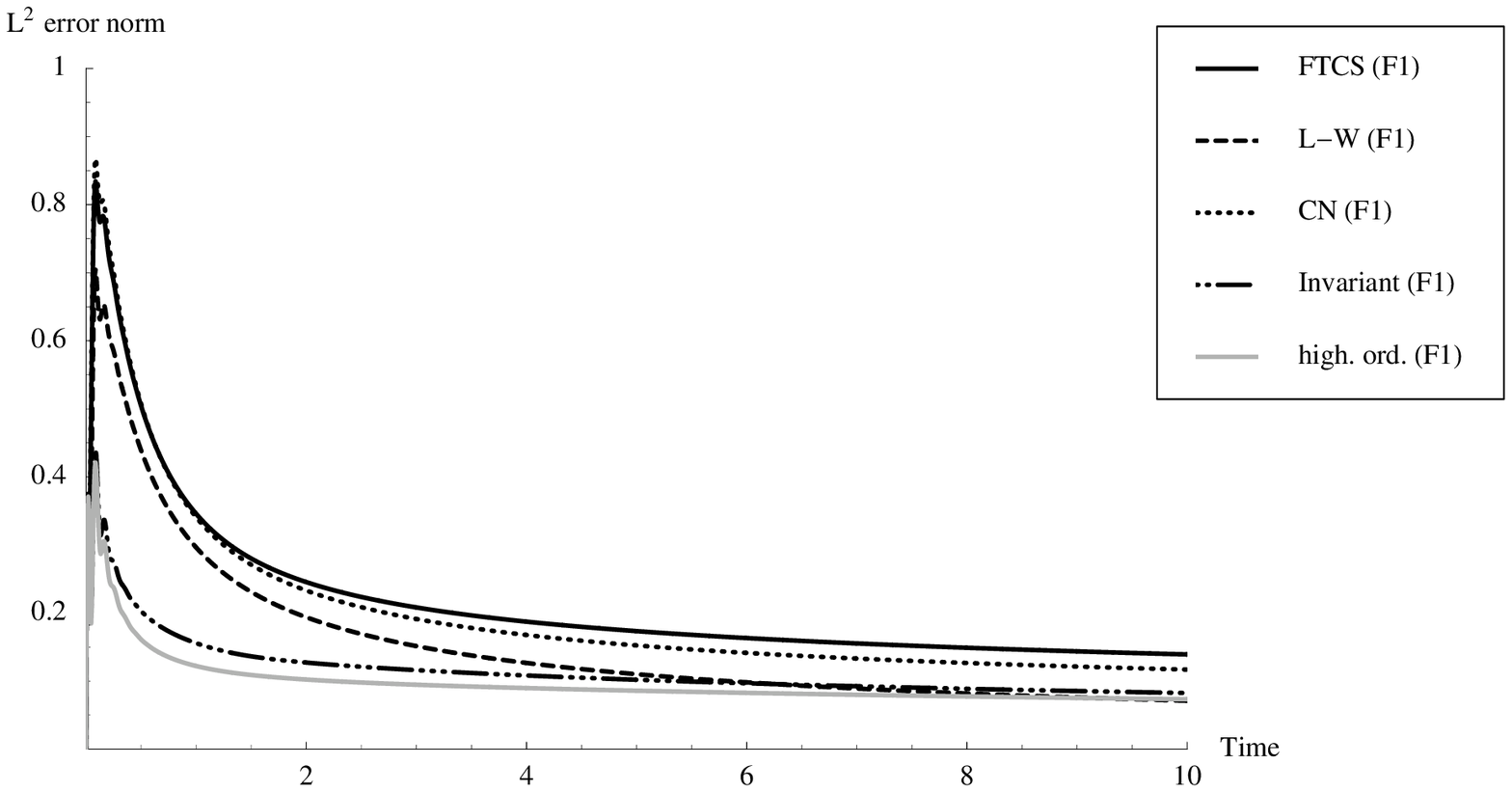}}
\resizebox*{0.49\columnwidth}{0.23\textheight}{\includegraphics{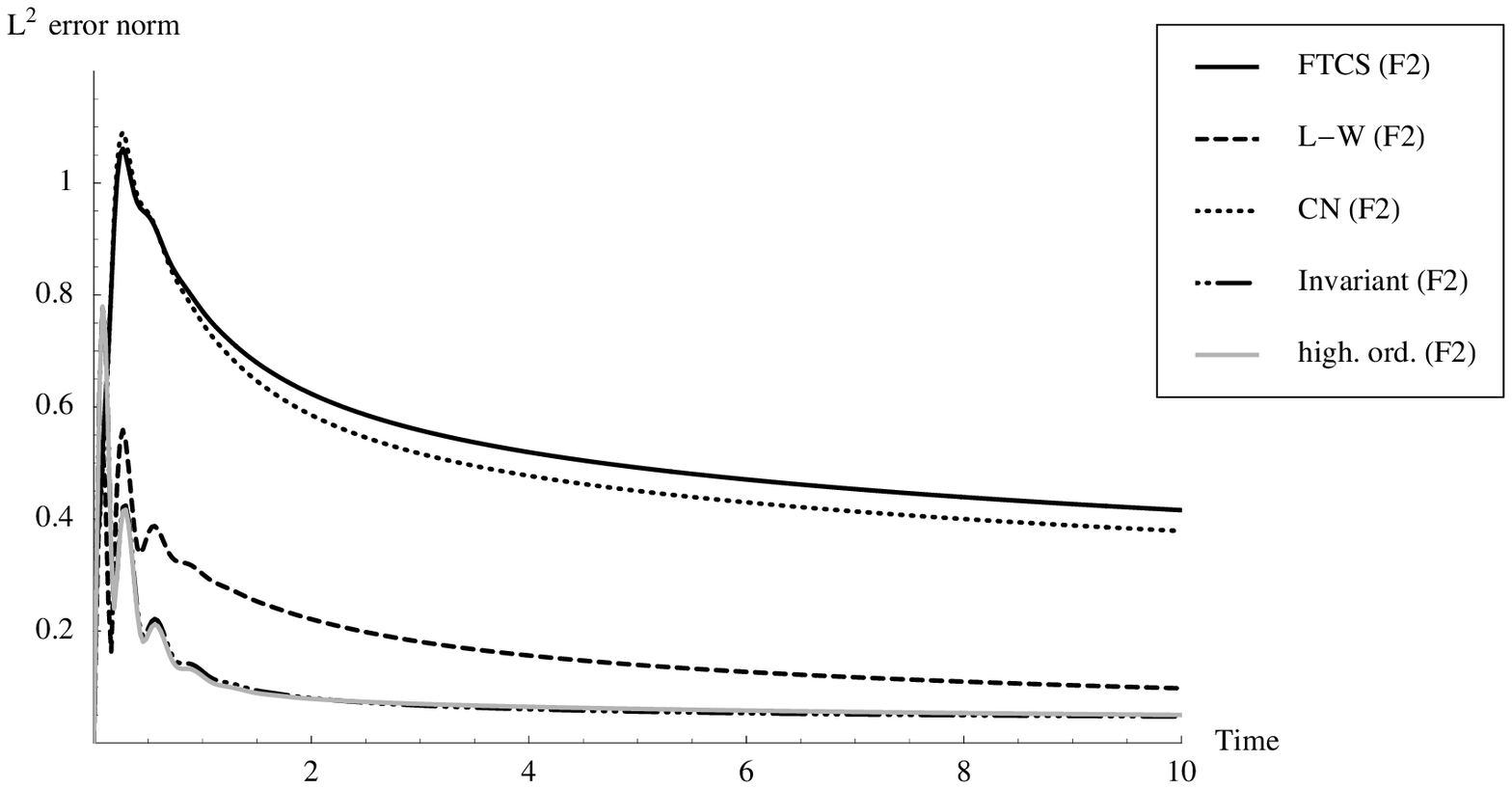}}
\caption{\footnotesize{Evolution of the error $L^2$-norm in (F1) and (F2). $Re_h=2$, $CFL=0.04$} \normalsize}
\label{L2NormTrG-2-004}
\end{figure}

\begin{figure}[!hbp]
\resizebox*{0.49\columnwidth}{0.23\textheight}{\includegraphics{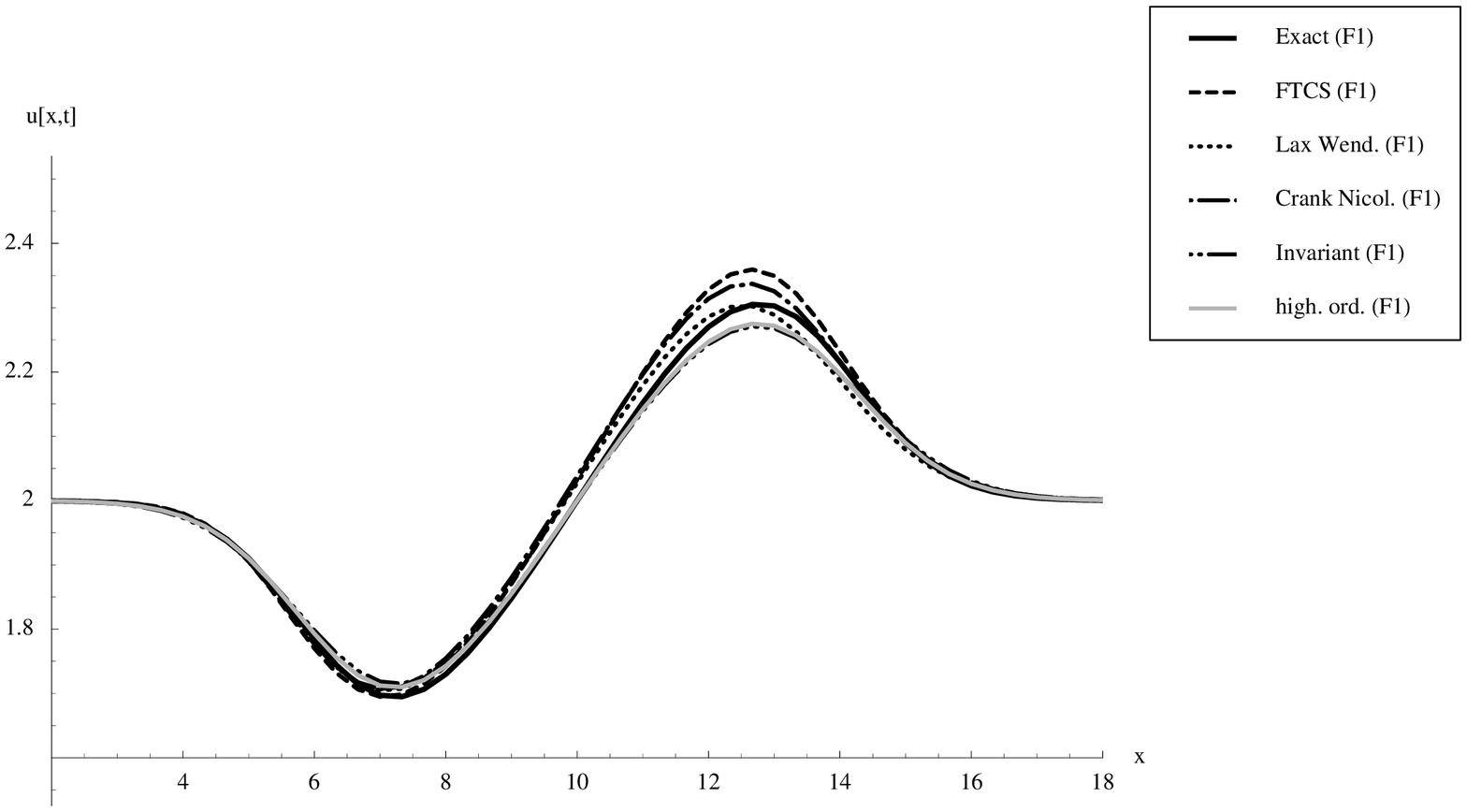}}
\resizebox*{0.49\columnwidth}{0.23\textheight}{\includegraphics{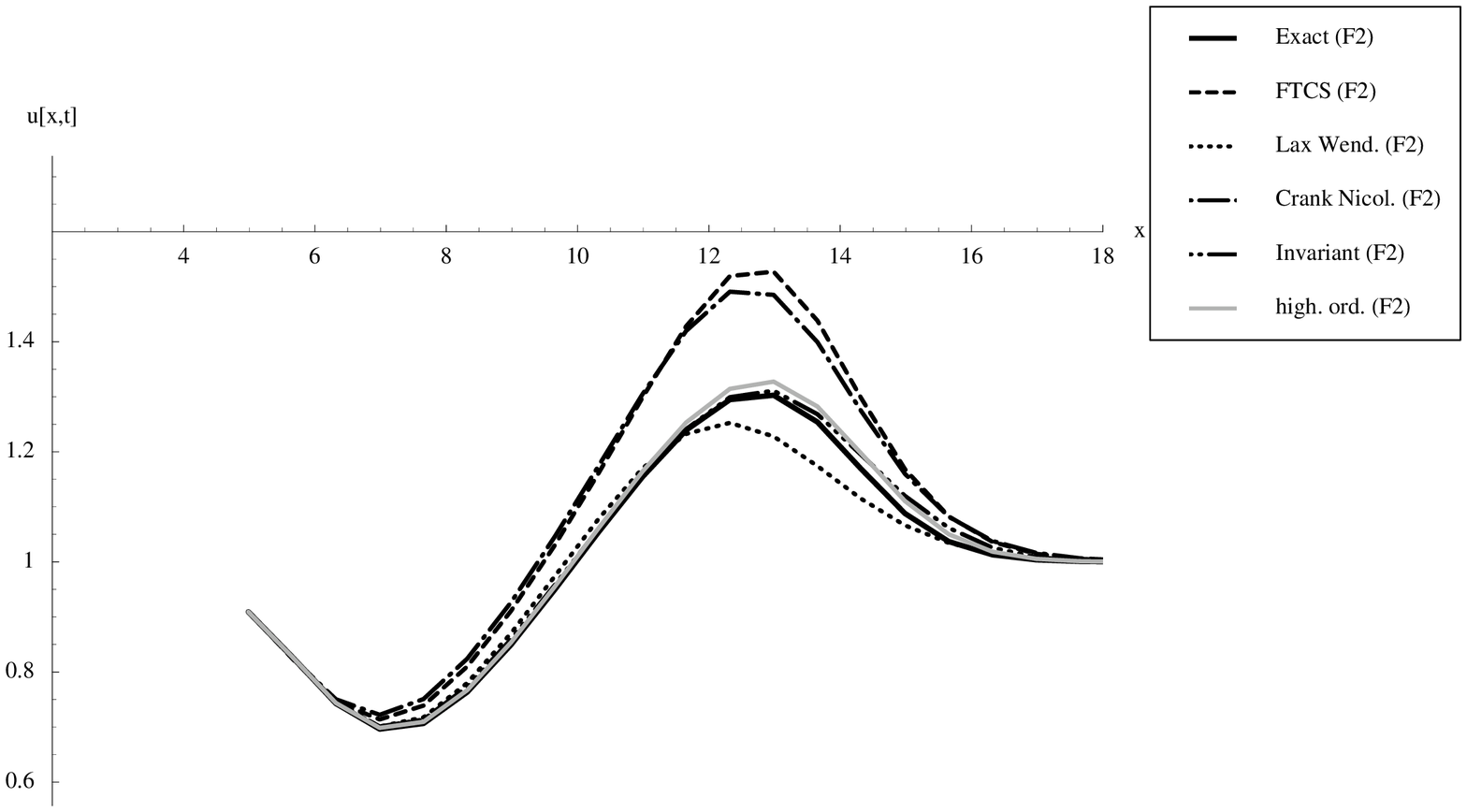}}
\caption{\footnotesize{Space variation of the numerical solutions of the schemes and the exact solution in $(F1)$ and $(F2)$.
$Re_h=2$, $CFL=0.04$} \normalsize}
\label{Scheme2-004}
\end{figure}

\begin{figure}[!tbp]
\resizebox*{0.49\columnwidth}{0.23\textheight}{\includegraphics{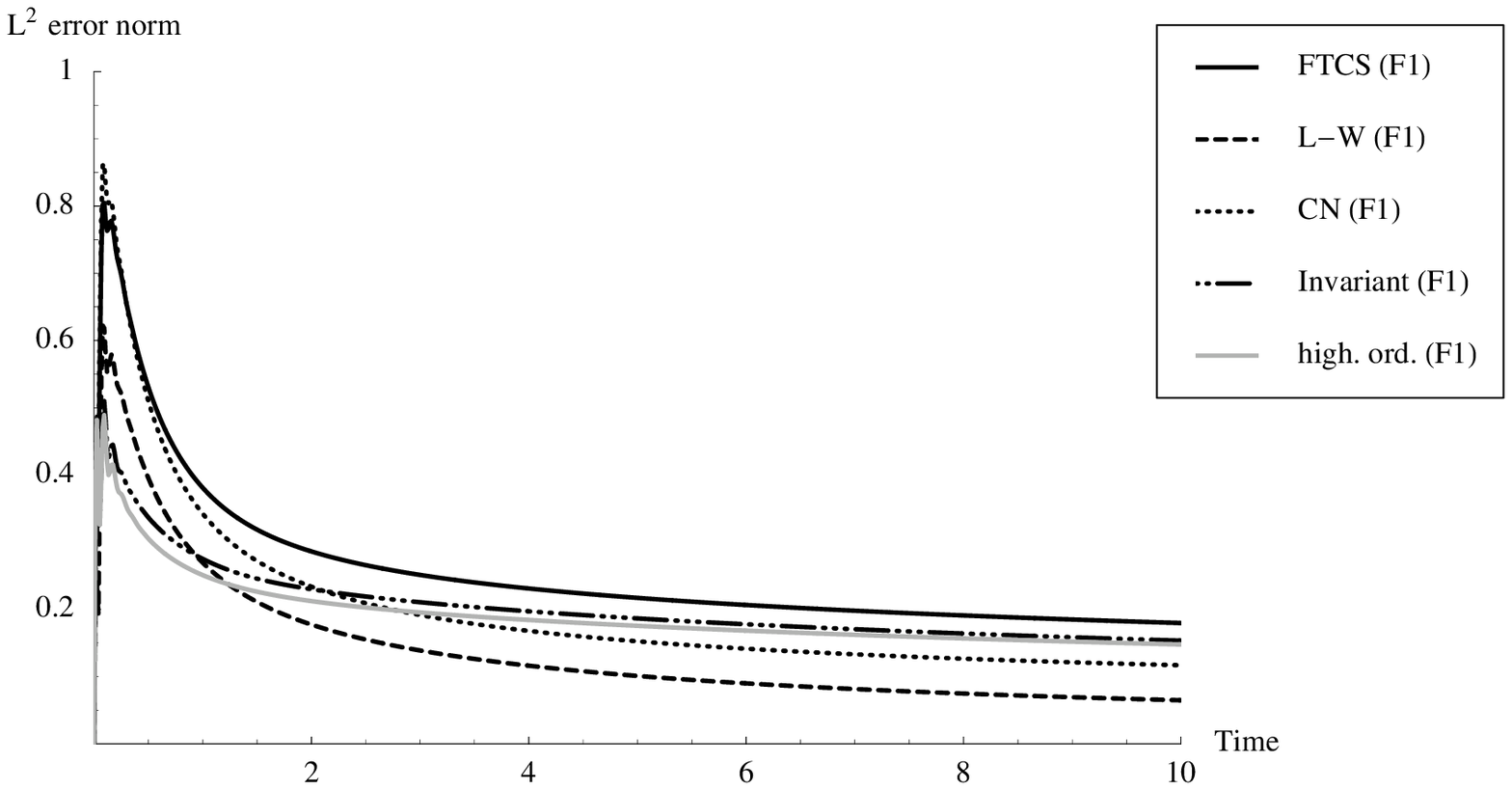}}
\resizebox*{0.49\columnwidth}{0.23\textheight}{\includegraphics{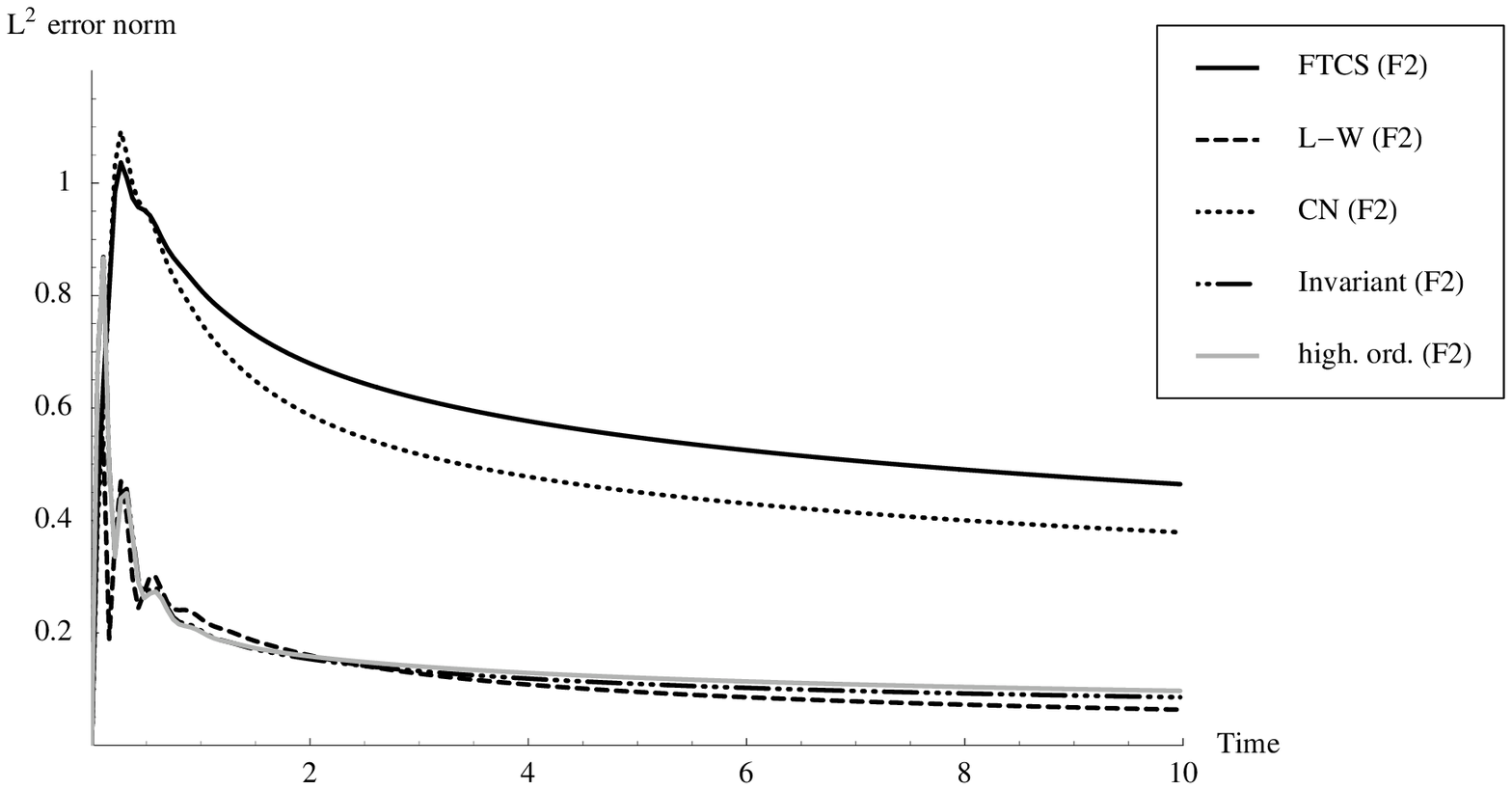}}
\caption{\footnotesize{Evolution of the error $L^2$-norm in (F1) and (F2). $Re_h=2$, $CFL=0.08$} \normalsize}
\label{L2NormTrG-2-008}
\end{figure}

\begin{figure}[!tbp]
\resizebox*{0.49\columnwidth}{0.23\textheight}{\includegraphics{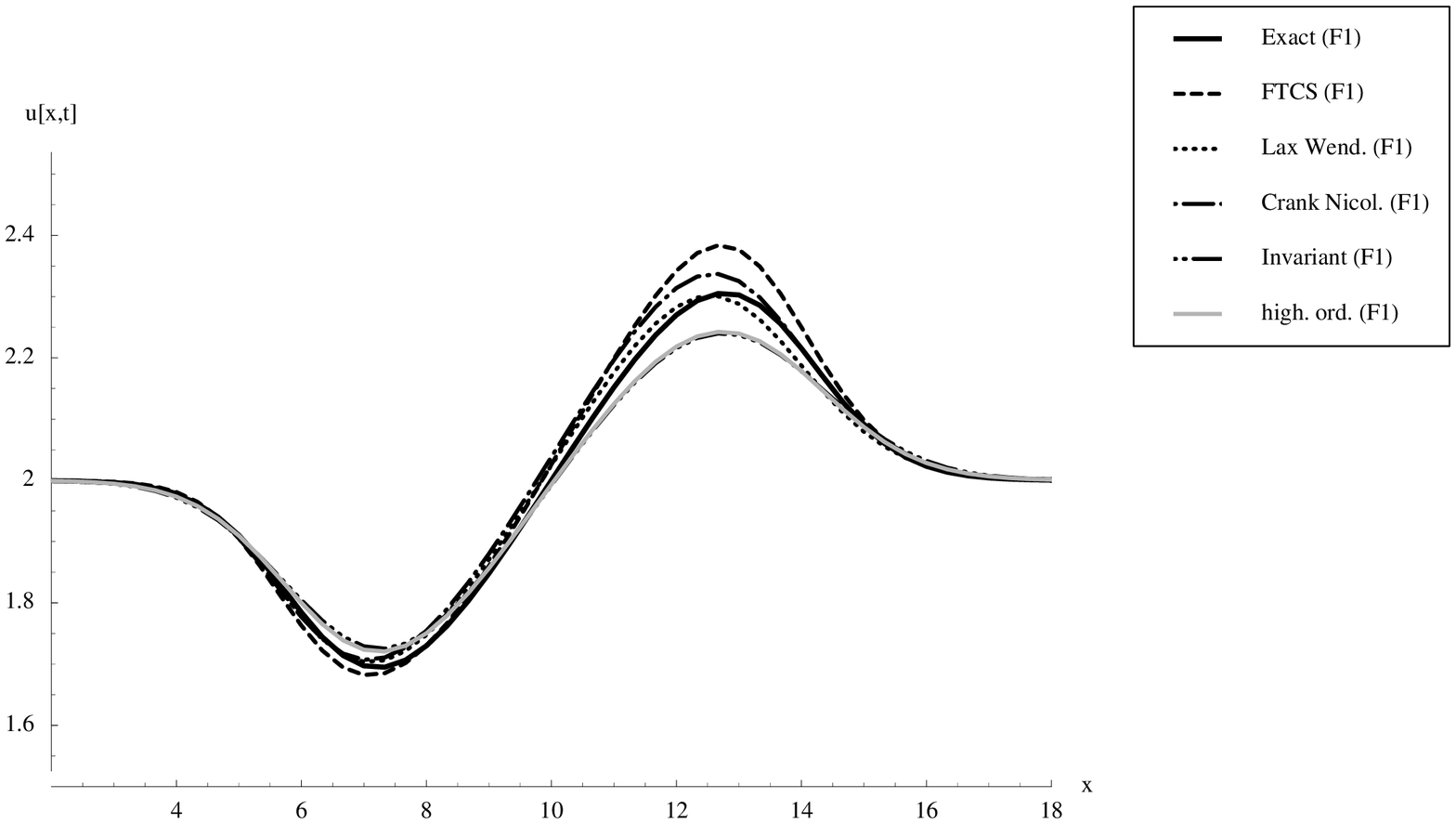}}
\resizebox*{0.49\columnwidth}{0.23\textheight}{\includegraphics{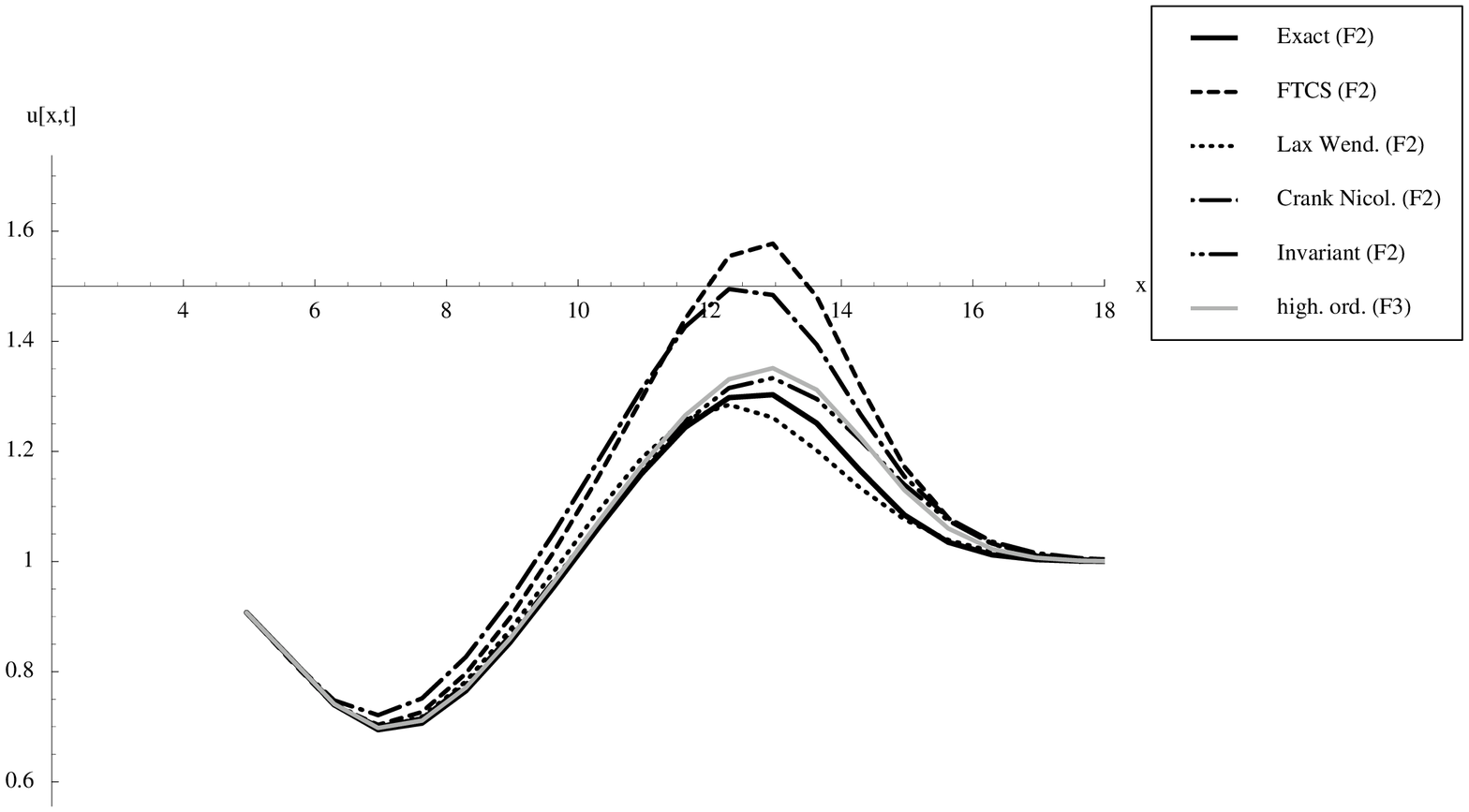}}
\caption{\footnotesize{Space variation of the numerical solutions of the schemes and the exact solution in $(F1)$ and $(F2)$.
$Re_h=2$, $CFL=0.08$} \normalsize}
\label{Scheme2-008}
\end{figure}

\begin{figure}[!tbp]
\resizebox*{0.49\columnwidth}{0.23\textheight}{\includegraphics{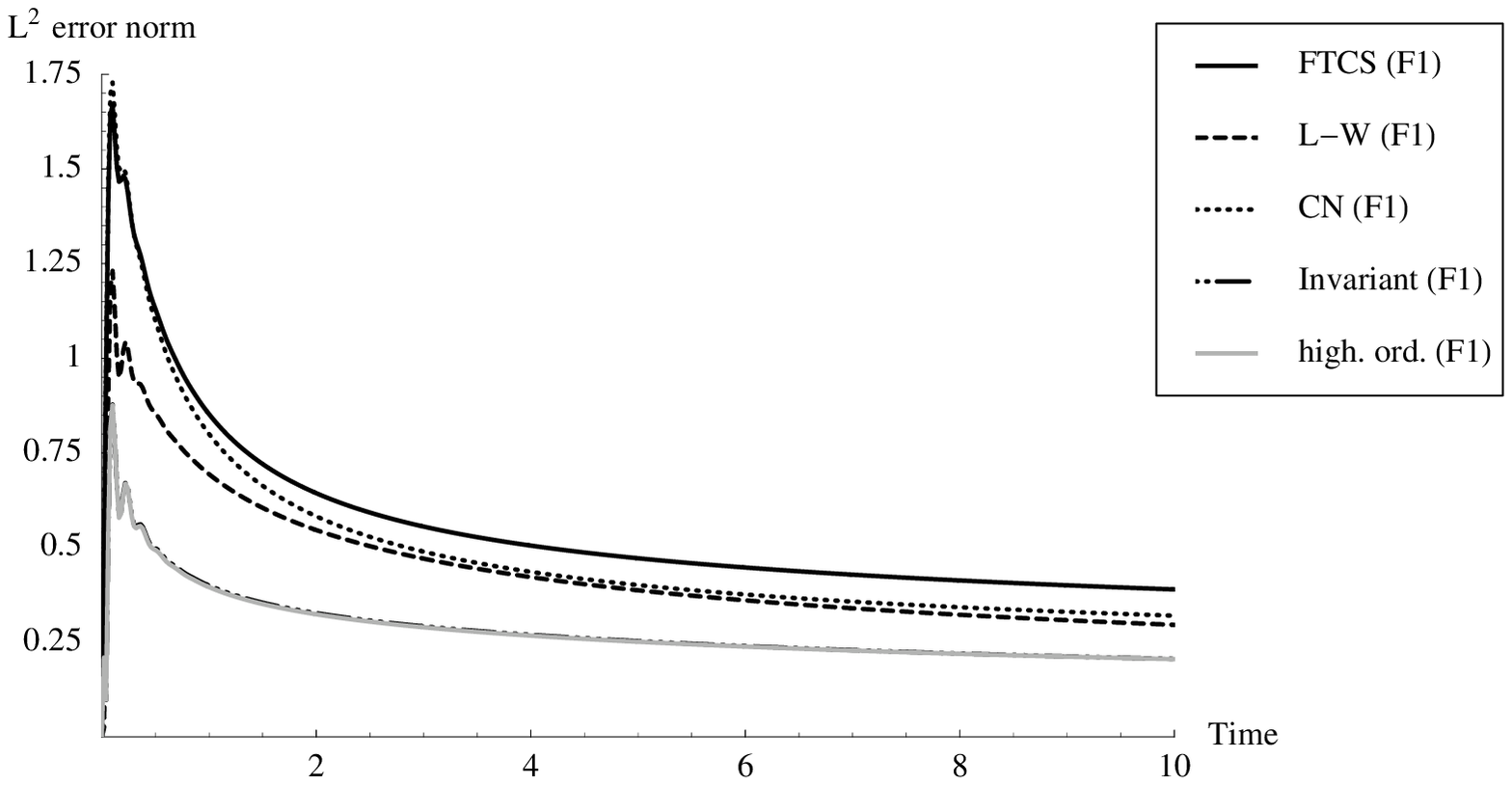}}
\resizebox*{0.49\columnwidth}{0.23\textheight}{\includegraphics{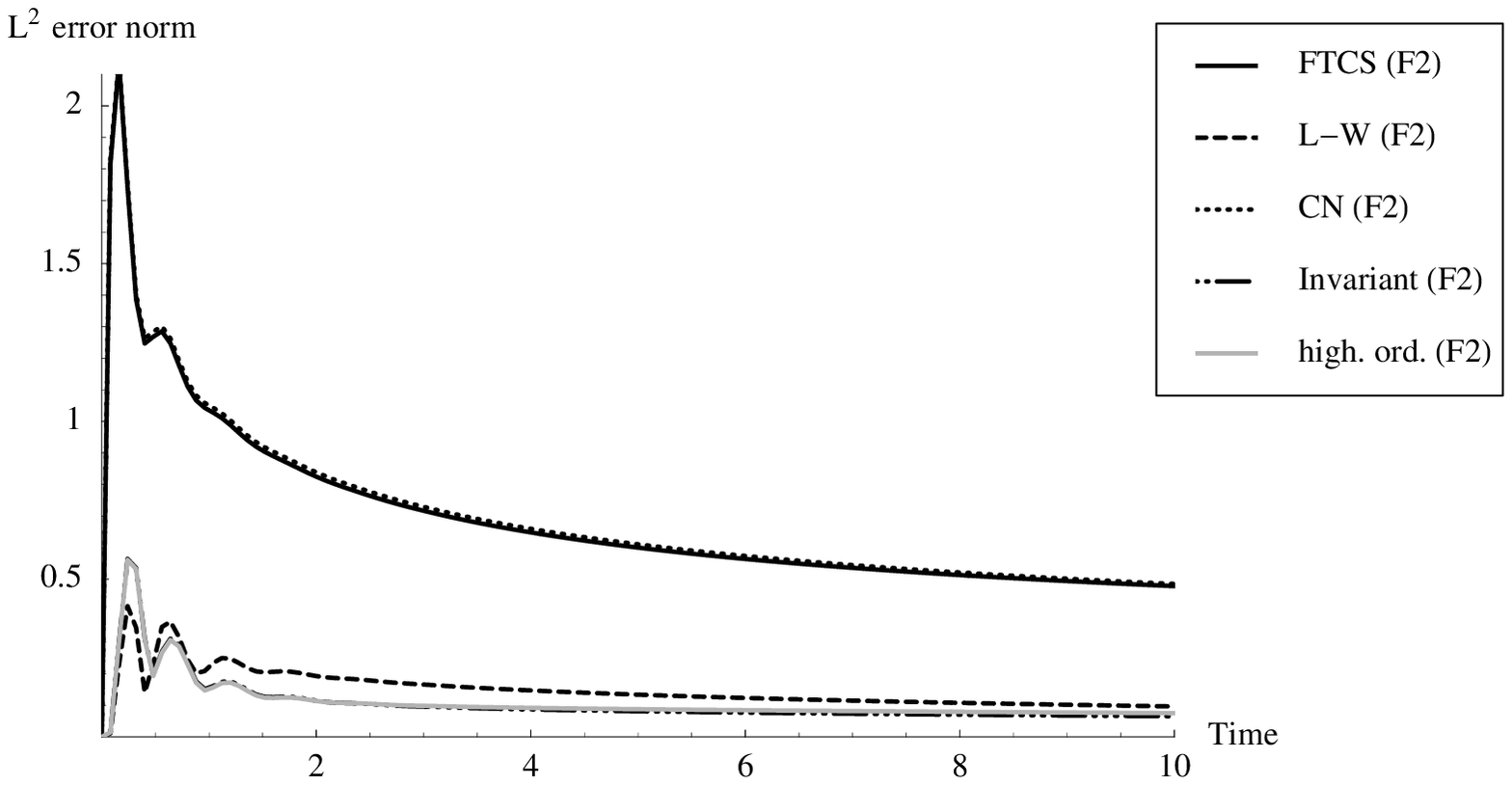}}
\caption{\footnotesize{Evolution of the error $L^2$-norm in (F1) and (F2). $Re_h=3$, $CFL=0.08$} \normalsize}
\label{L2NormTrG-3-008}
\end{figure}

\begin{figure}[!hbp]
\resizebox*{0.49\columnwidth}{0.23\textheight}{\includegraphics{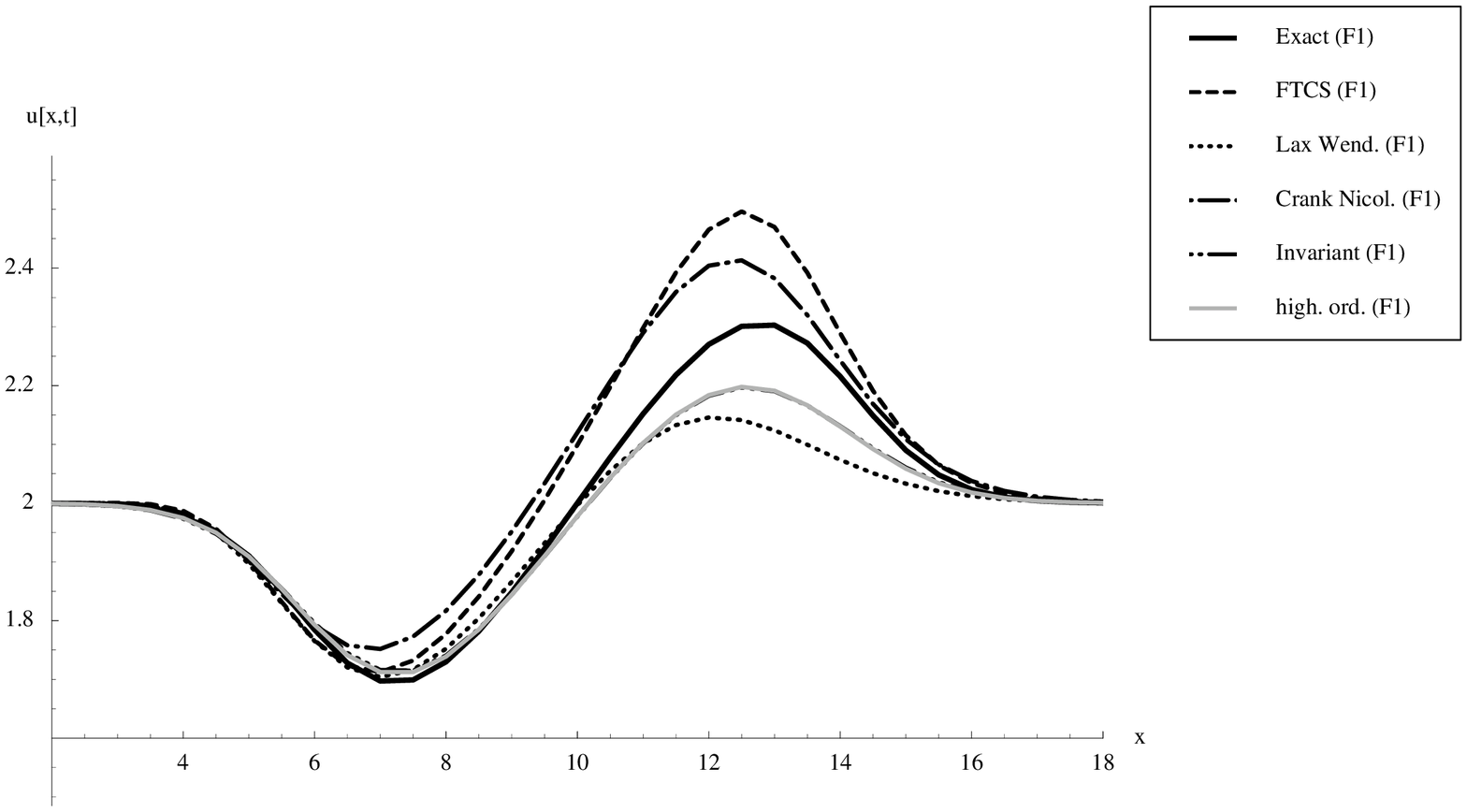}}
\resizebox*{0.49\columnwidth}{0.23\textheight}{\includegraphics{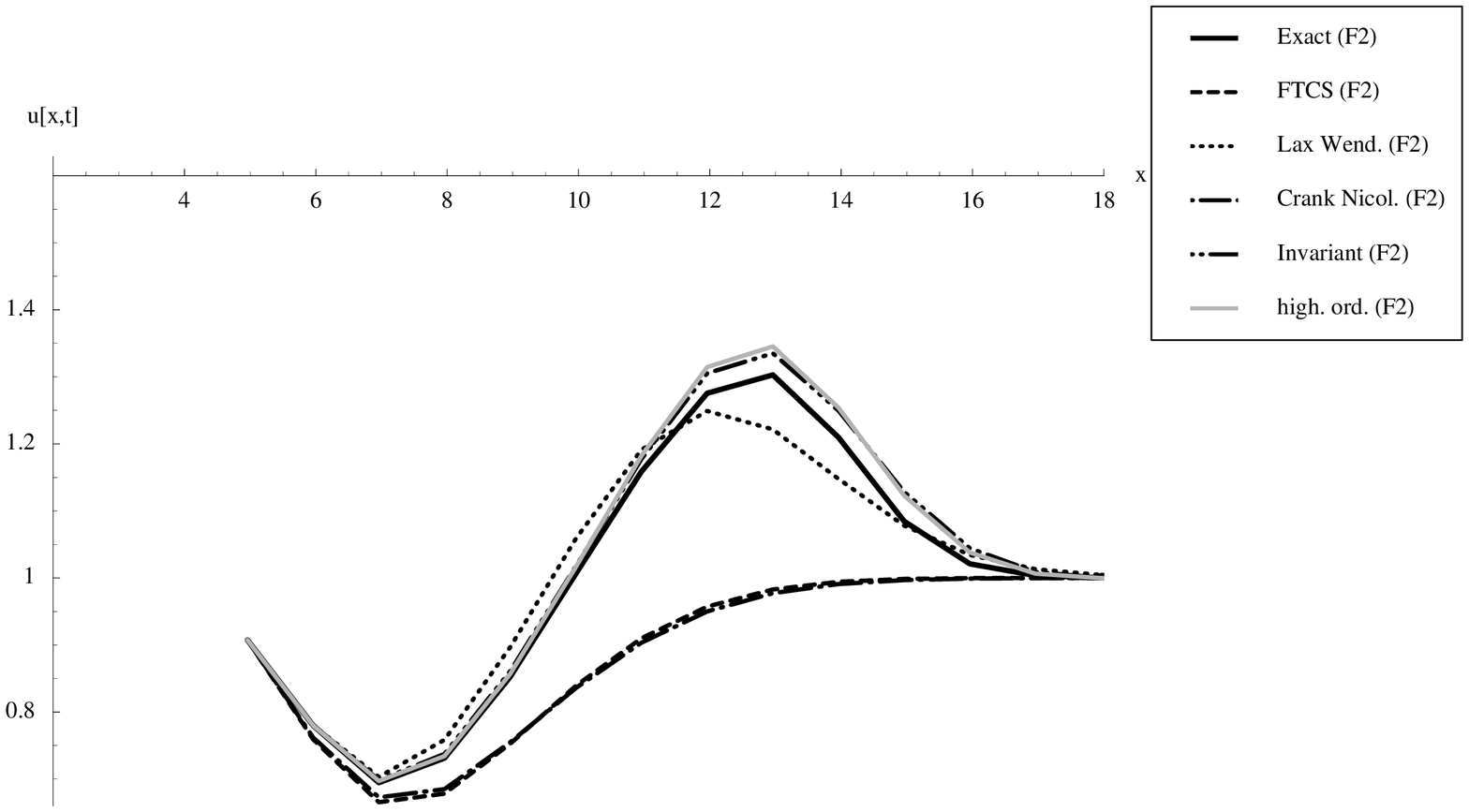}}
\caption{\footnotesize{Space variation of the schemes numerical and
of the exact solution in $(F1)$ and $(F2)$. $Re_h=3$, $CFL=0.08$}
\normalsize} \label{Scheme3-008}
\end{figure}

\noindent The problem consists in solving the following differential system:
\footnotesize{
\begin{eqnarray}
\displaystyle & &u_t+u u_x-\nu u_{xx}=0,\ x\in[0,40],\ t\in[0,20]\nonumber\\
\displaystyle & &u(x,0)=f(x),\ \displaystyle u(0,t)=g(t),\
u(40,t)=h(t)\nonumber
\end{eqnarray}
} \normalsize \noindent The initial and boundary conditions, $f$, $h$, and $g$
are provided by an exact solution of the Burgers equation:
\footnotesize{
\begin{eqnarray}
\displaystyle u(x,t)=\frac{(x-2 t)/(t+0.1)}{1+\nu^2\sqrt{t+0.1}\exp{((x-2t)^2/(4 \nu(t+0.1)))}}+2
\end{eqnarray}
} \normalsize \noindent Figures \ref{L2NormTrG-2-004},
\ref{L2NormTrG-2-008} and \ref{L2NormTrG-3-008} show the time
evolution of the $L^2$-norm of the error for the considered schemes, for specific
values of the $CFL$ number and the mesh Reynolds number $Re_h$.
Figures \ref{Scheme2-004},  \ref{Scheme2-008} and \ref{Scheme3-008} 
display the variations, as functions of the space variable, of the numerical 
solutions of the considered schemes for the specific value $t=5$.
In each frame, the numerical solutions are compared to the exact one.

\noindent The error analysis of the invariant scheme in the reference frame 
through the features of the truncation error and the graphical representation of 
the norms of the error (cf. Figures \ref{L2NormTrG-2-004}, \ref{L2NormTrG-2-008} and \ref{L2NormTrG-3-008}) 
allows to say that the invariant scheme is dissipative and slightly dispersive.

\noindent The presence of the dissipative term $(C u_x)_x$ in the differential representation 
of the invariant scheme and the presence of the higher order error terms involving 
the even-order derivative $u_{6x}$ (cf. Equation (\ref{eqn:4-order-acc})) show that the
scheme produces numerical damping. Particularly, the amplitudes are not correctly 
represented for high frequencies, since the solution is subjected to rather rough 
variation during the first iterations. The dissipation is stronger for $Re_h=2$, $CFL=0.08$ 
in the reference frame (see Figure \ref{Scheme2-008}). Moreover, the presence of higher order
error terms involving the odd-order derivative $u_{5x}$ corresponds to a phase error.

\noindent The non-invariant schemes are more altered by the change of the frame than the invariant one.
Moreover, the invariant scheme appears to be as accurate as the higher order one in the frame (F2).


\newpage

\medskip


\end{document}